\documentclass{amsart}
\usepackage{amsmath,amssymb}
\usepackage[dvips]{graphics}
\usepackage[all]{xy}
\newtheorem{Theorem}{Theorem}[section]
\newtheorem{Lemma}[Theorem]{Lemma}
\newtheorem{Proposition}[Theorem]{Proposition}

\newtheorem{Remark}[Theorem]{Remark}

\makeatletter
\@addtoreset{figure}{section}%{subsection}
\def\@thmcountersep{-}
\makeatother

\numberwithin{equation}{section}

\begin{document}

\title{Circle immersions that can be divided into two arc embeddings}

\author{Kouki Taniyama}
\address{Department of Mathematics, School of Education, Waseda University, Nishi-Waseda 1-6-1, Shinjuku-ku, Tokyo, 169-8050, Japan}
%\curraddr{}
\email{taniyama@waseda.jp}
\thanks{The third author was partially supported by Grant-in-Aid for Scientific Research (C) (No. 18540101), Japan Society for the Promotion of Science.}

%    General info

\date{}

\dedicatory{}


\begin{abstract}
We give a complete characterization of a circle immersion that can be divided into two arc embeddings in terms of its chord diagram.
\end{abstract}

\maketitle

\section{Introduction} 

Let ${\mathbb S}^1$ be the unit circle. Let $X$ be a set and $f:{\mathbb S}^1\to X$ a map. Let $n$ be a natural number greater than one. Suppose that there are $n$ subspaces $I_1,\cdots,I_n$ of ${\mathbb S}^1$ with the following properties.

(1) Each $I_i$ is homeomorphic to a closed interval.

(2) ${\mathbb S}^1=I_1\cup\cdots\cup I_n$.

(3) The restriction map $f|_{I_i}:I_i\to X$ is injective for each $i$.

\noindent
Then we say that $f$ can be divided into $n$ arc embeddings. We define the {\it arc number} of $f$, denoted by ${\rm arc}(f)$, to be the smallest such $n$ except the case that $f$ itself is injective. If $f$ itself is injective then we define ${\rm arc}(f)=1$. If $f$ cannot be divided into $n$ arc embeddings for any natural number $n$ then we define ${\rm arc}(f)=\infty$. 
Note that if $f$ can be divided into $n$ arc embeddings then there exist $n$ subspaces $I_1,\cdots,I_n$ of ${\mathbb S}^1$ with (1), (2) and (3) above together with the following additional condition.

(4) $I_i\cap I_j=\partial I_i\cap\partial I_j$ for each $i$ and $j$ with $1\leq i<j\leq n$.

\noindent
Namely we may assume that ${\mathbb S}^1$ is covered by mutually interior disjoint $n$ simple arcs $I_1,\cdots,I_n$.

Let ${\rm S}(f)=\{x\in{\mathbb S}^1|f^{-1}(f(x)){\rm\ is\ not\ a\ singleton.}\}$ and ${\rm s}(f)=f({\rm S}(f))$.
We say that a map $f:{\mathbb S}^1\to X$ has {\it finite multiplicity} if ${\rm S}(f)$ is a finite subset of ${\mathbb S}^1$.
From now on we restrict our attention to maps that have finite multiplicity. The purpose of this paper is to give a characterization of a map $f:{\mathbb S}^1\to X$ with ${\rm arc}(f)=2$. By $|Y|$ we denote the cardinality of a set $Y$. Let $m(f)$ be the maximum of $|f^{-1}(y)|$ where $y$ varies over all points of $X$.
It is clear that ${\rm arc}(f)\geq m(f)$. 
Thus we further restrict our attention to a map $f:{\mathbb S}^1\to X$ whose multiple points are only finitely many double points. 
Namely $f$ has finite multiplicity and $m(f)\leq2$. Then we have $|{\rm S}(f)|=2m$ for some non-negative integer $m$. Then the {\it crossing number} of $f$, denoted by $c(f)$, is defined by $m$. 

\begin{figure}[htbp]
      \begin{center}
\scalebox{0.4}{\includegraphics*{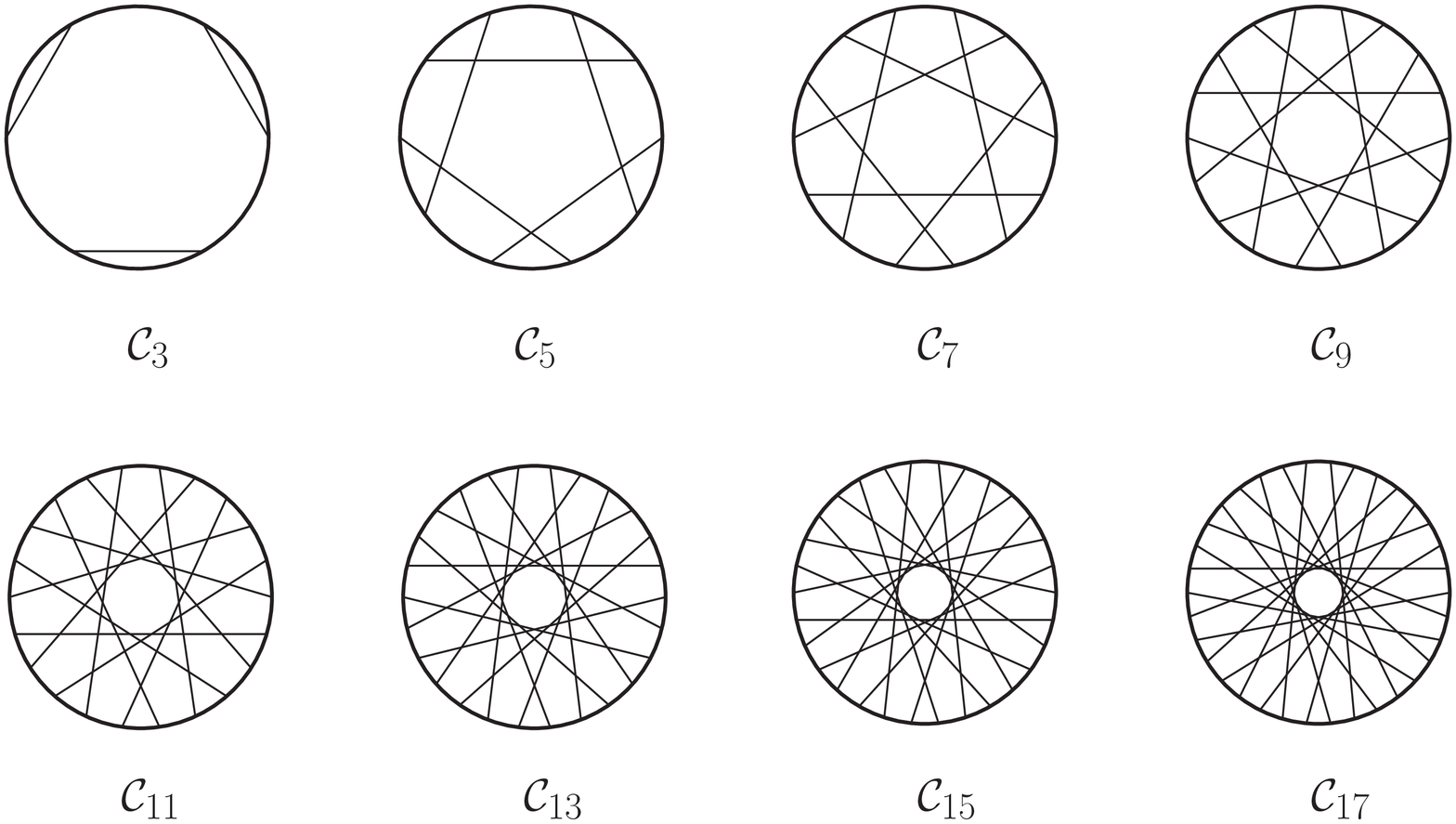}}
      \end{center}
   \caption{}
  \label{chord-diagrams}
\end{figure}

Let $m$ be a natural number. An {\it $m$-chord diagram} on ${\mathbb S}^1$ is a pair ${\mathcal C}=(P,\varphi)$ where $P$ is a subset of ${\mathbb S}^1$ that contains exactly $2m$ points and $\varphi$ is a fixed point free involution on $P$. A {\it chord} $c$ of ${\mathcal C}$ is an unordered pair of points $(x,\varphi(x))=(\varphi(x),x)$ where $x$ is a point in $P$. 
Let $\sim_{\mathcal C}$ be the equivalence relation on ${\mathbb S}^1$ generated by $x\sim_{\mathcal C}\varphi(x)$ for every $x\in P$. Let ${\mathbb S}^1/\sim_{\mathcal C}$ be the quotient space and $f_{\mathcal C}:{\mathbb S}^1\to {\mathbb S}^1/\sim_{\mathcal C}$ the quotient map. We call $f_{\mathcal C}$ the associated map of ${\mathcal C}$. Then the arc number of ${\mathcal C}$, denoted by ${\rm arc}({\mathcal C})$, is defined to be the arc number of $f_{\mathcal C}$. 
Two $m$-chord diagrams ${\mathcal C}_1=(P_1,\varphi_1)$ and ${\mathcal C}_2=(P_2,\varphi_2)$ are {\it equivalent} if there is an orientation preserving self-homeomorphism $h$ of ${\mathbb S}^1$ such that $h(P_1)=P_2$ and $h\circ\varphi_1=\varphi_2\circ h$. From now on we consider $m$-chord diagrams up to this equivalence relation. 
In the following we sometimes express an $m$-chord diagram ${\mathcal C}=(P,\varphi)$ by $m$ line-segments in the plane ${\mathbb R}^2$ where ${\mathbb S}^1\subset{\mathbb R}^2$, $P$ is the set of the end points of these line-segments and $x$ and $\varphi(x)$ are joined by a line segment for each $x\in P$. Thus a line segment express a chord and from now on we do not distinguish them. 
See for example Figure \ref{chord-diagrams}.

Let $f:{\mathbb S}^1\to X$ be a map whose multiple points are only finitely many double points. By ${\mathcal C}(f)$ we denote the $c(f)$-chord diagram $({\rm S}(f),\varphi_f)$ where $\varphi_f:{\rm S}(f)\to {\rm S}(f)$ is the fixed point free involution with $f|_{{\rm S}(f)}\circ\varphi_f=f|_{{\rm S}(f)}$. We call ${\mathcal C}(f)$ the associated chord diagram of $f$. Then it is clear that ${\rm arc}(f)={\rm arc}({\mathcal C}(f))$. 

A chord diagram ${\mathcal D}=(Q,\psi)$ is called a {\it sub-chord diagram} of a chord diagram ${\mathcal C}=(P,\varphi)$ if $Q$ is a subset of $P$ and $\psi$ is the restriction of $\varphi$ on $Q$. Then it is clear that ${\rm arc}({\mathcal D})\leq{\rm arc}({\mathcal C})$. We call ${\mathcal D}=(Q,\psi)$ a {\it proper sub-chord diagram} of ${\mathcal C}=(P,\varphi)$ if ${\mathcal D}$ is a sub-chord diagram of ${\mathcal C}$ and $Q$ is a proper subset of $P$. 

Let $n$ be a natural number. Let ${\mathcal C}_{2n+1}$ be a $(2n+1)$-chord diagram as illustrated in Figure \ref{chord-diagrams}. 
To give a more precise definition we introduce the followings. 
Let $k$ be a natural number greater than two. Let $R_k$ be a regular $k$-gon inscribed in ${\mathbb S}^1$ and $v_{k;1},\cdots,v_{k;k}$ the vertices of $R_k$ that are arranged in this order on ${\mathbb S}^1$ along the counterclockwise orientation of ${\mathbb S}^1$. 
Namely $v_{k;i}$ and $v_{k;i+1}$ are adjacent in $R_k$ for each $i$ where the indices are considered modulo $k$. Let $j$ is a natural number less than $\frac{k}{2}$. Let $c(k;i,j)$ be the chord joining $v_{k;i}$ and $v_{k;i+j}$ for each $i\in\{1,\cdots,k\}$. 
Then ${\mathcal C}_{2n+1}$ is the chord diagram represented by chords $c(4n+2;2i-1,2n-1)$ with $i\in\{1,\cdots,2n+1\}$. 
We will show that ${\rm arc}({\mathcal C}_{2n+1})=3$ but ${\rm arc}({\mathcal D})=2$ for any proper sub-chord diagram ${\mathcal D}$ of ${\mathcal C}_{2n+1}$. 
Then we have the following theorem.

\vskip 3mm

\begin{Theorem}\label{main-theorem}
Let $m$ be a natural number and ${\mathcal C}$ an $m$-chord diagram on ${\mathbb S}^1$. Then ${\rm arc}({\mathcal C})=2$ if and only if no sub-chord diagram of ${\mathcal C}$ is equivalent to the chord diagram ${\mathcal C}_{2n+1}$ for any natural number $n$.
\end{Theorem}

\vskip 3mm

The motive for this paper was the result in \cite{Hotz} that every knot has a diagram which can be divided into two simple arcs. This result is re-discovered by \cite{Ozawa} and \cite{Shinjo}. See also \cite{A-S-T}. Then it is natural to ask what plane closed curve can be divided into two simple arcs. Theorem \ref{main-theorem} gives an answer to this question. However we still have a question whether or not do we actually need all of ${\mathcal C}_3,{\mathcal C}_5,\cdots$. 
The following proposition answers this question that we actually need all of them.

\vskip 3mm

\begin{Proposition}\label{example}
For each natural number $n$ there exist a smooth immersion $f_n:{\mathbb S}^1\to {\mathbb R}^2$ with ${\rm arc}(f_n)=3$ that has only finitely many transversal double points such that the associated chord diagram ${\mathcal C}(f_n)$ of $f_n$ has a sub-chord diagram which is equivalent to ${\mathcal C}_{2n+1}$ but has no sub-chord diagram which is equivalent to ${\mathcal C}_{2m+1}$ for any $m<n$. 
\end{Proposition}

\vskip 3mm

%\vskip 3mm
% 

%
\section{Proof of Theorem 1.1} 

First we check that ${\rm arc}({\mathcal C}_{2n+1})=3$. Let ${\mathcal C}=(P,\varphi)$ be an $m$-chord diagram with ${\rm arc}({\mathcal C})=2$. A pair of points $p,q\in{\mathbb S}^1\setminus P$ are called a {\it cutting pair} of ${\mathcal C}$ if $x$ and $\varphi(x)$ belong to the different components of ${\mathbb S}^1\setminus\{p,q\}$ for each $x\in P$. Then we have that $p$ and $q$ are \lq\lq antipodal\rq\rq. Namely we have that each component of ${\mathbb S}^1\setminus\{p,q\}$ contains exactly $m$ points of $P$. Thus we can check whether or not a given $m$-chord diagram has arc number $2$ by examining $m$ pairs of antipodal points of it. Then by the symmetry of ${\mathcal C}_{2n+1}$ we immediately have that ${\rm arc}({\mathcal C}_{2n+1})>2$. Then it is easily seen that ${\rm arc}({\mathcal C}_{2n+1})=3$ and ${\rm arc}({\mathcal D})=2$ for any proper sub-chord diagram ${\mathcal D}$ of ${\mathcal C}_{2n+1}$. Then the \lq only if part\rq\  of the proof of Theorem \ref{main-theorem} immediately follows. 
The \lq if part\rq\ immediately follows from the following proposition.

\vskip 3mm

\begin{Proposition}\label{proposition}
Let ${\mathcal C}$ be a chord diagram on ${\mathbb S}^1$ that satisfies the following condition $(\ast)$.

$(\ast)$ ${\rm arc}({\mathcal C})\geq3$ and ${\rm arc}({\mathcal D})=2$ for any proper sub-chord diagram ${\mathcal D}$ of ${\mathcal C}$. 

\noindent
Then there is a natural number $n$ such that ${\mathcal C}$ is equivalent to ${\mathcal C}_{2n+1}$.
\end{Proposition}

\vskip 3mm

Note that deleting a chord will decrease the arc number at most by one. Therefore, if ${\mathcal C}$ is a chord diagram on ${\mathbb S}^1$ that satisfies the condition $(\ast)$, then ${\rm arc}({\mathcal C})=3$.
For the proof of Proposition \ref{proposition} we prepare the following lemmas. 
Let ${\mathcal C}=(P,\varphi)$ be a chord diagram and $c=(x,\varphi(x))$ a chord of ${\mathcal C}$. Let $\alpha$ and $\beta$ be the components of ${\mathbb S}^1\setminus \{x,\varphi(x)\}$. We may suppose without loss of generality that $|\alpha\cap P|\leq|\beta\cap P|$. Then the {\it length} of $c$ in ${\mathcal C}$, denoted by $l(c)=l(c,{\mathcal C})$, is defined to be $|\alpha\cap P|+1$. By ${\mathcal C}\setminus c$, we denote the chord diagram $(P\setminus\{x,\varphi(x)\},\varphi|_{P\setminus\{x,\varphi(x)\}})$. Let $p,q,x$ and $y$ be mutually distinct four points on ${\mathbb S}^1$. We say that the pair of points $p$ and $q$ separates the pair of points $x$ and $y$ if each component of ${\mathbb S}^1\setminus \{p,q\}$ contains exactly one of $x$ and $y$. Note that the pair of points $p$ and $q$ separates the pair of points $x$ and $y$ if and only if the chord joining $p$ and $q$ intersects the chord joining $x$ and $y$.

\vskip 3mm

\begin{Lemma}\label{lemma1}
Let ${\mathcal C}=(P,\varphi)$ be a chord diagram on ${\mathbb S}^1$ that satisfies the condition $(\ast)$. Let $c=(x,\varphi(x))$ be a chord of ${\mathcal C}$. Let $p$ and $q$ be a cutting pair of ${\mathcal C}\setminus c$. Then $p$ and $q$ do not separate $x$ and $\varphi(x)$.
\end{Lemma}

\vskip 3mm

\noindent{\bf Proof.}
If $p$ and $q$ separate $x$ and $\varphi(x)$, then $p$ and $q$ is a cutting pair of ${\mathcal C}$ itself. Then it follows ${\rm arc}({\mathcal C})=2$. This is a contradiction. $\Box$

\vskip 3mm

\begin{Lemma}\label{lemma2}
Let ${\mathcal C}=(P,\varphi)$ be an $m$-chord diagram on ${\mathbb S}^1$ that satisfies the condition $(\ast)$. 
Let $c=(x,\varphi(x))$ be a chord of ${\mathcal C}$. Then $l(c,{\mathcal C})\leq m-2$.
\end{Lemma}

\vskip 3mm

\noindent{\bf Proof.}
Since $|P|=2m$ we have $1\leq l(c,{\mathcal C})\leq m$ for any chord $c$. First we examine the case $l(c,{\mathcal C})=m$. In this case we have that each component of ${\mathbb S}^1\setminus \{x,\varphi(x)\}$ contains exactly $m-1$ elements of $P$. Let $p$ and $q$ be a cutting pair of ${\mathcal C}\setminus c$. Then by Lemma \ref{lemma1} we have that $p$ and $q$ do not separate $x$ and $\varphi(x)$. Note that each component of ${\mathbb S}^1\setminus \{p,q\}$ also contains exactly $m-1$ elements of $P$. Then it follows that $p$ and $q$ are next to $x$ and $\varphi(x)$ or $\varphi(x)$ and $x$ respectively. We may suppose without loss of generality that $p$ and $q$ are next to $x$ and $\varphi(x)$ respectively.
Let $p'$ be a point on ${\mathbb S}^1$ that is next to $x$ and such that $p'$ and $p$ separate $x$ and $\varphi(x)$. Then we have that $p'$ and $q$ is a cutting pair of ${\mathcal C}$. This is a contradiction. 
Next we examine the case $l(c,{\mathcal C})=m-1$. In this case we have that one component of ${\mathbb S}^1\setminus \{x,\varphi(x)\}$ contains exactly $m-2$ elements of $P$ and the other component contains exactly $m$ elements of $P$. Let $p$ and $q$ be a cutting pair of ${\mathcal C}\setminus c$. Then by Lemma \ref{lemma1} we have that $p$ and $q$ do not separate $x$ and $\varphi(x)$. Note that each component of ${\mathbb S}^1\setminus \{p,q\}$  contains exactly $m-1$ elements of $P$. Then it follows that one of $p$ and $q$, say $p$ is next to $x$ or $\varphi(x)$, say $x$. Let $p'$ be a point on ${\mathbb S}^1$ that is next to $x$ and such that $p'$ and $p$ separate $x$ and $\varphi(x)$. Then we have that $p'$ and $q$ is a cutting pair of ${\mathcal C}$. This is a contradiction. 
Thus we have $l(c,{\mathcal C})\leq m-2$. $\Box$

\vskip 3mm

\begin{Lemma}\label{lemma3}
Let ${\mathcal C}=(P,\varphi)$ be an $m$-chord diagram on ${\mathbb S}^1$ that satisfies the condition $(\ast)$. 
Let $c=(x,\varphi(x))$ be a chord of ${\mathcal C}$. Then $l(c,{\mathcal C})\geq m-2$.
\end{Lemma}

\vskip 3mm

\noindent{\bf Proof.}
Suppose that there is a chord $c=(x,\varphi(x))$ of ${\mathcal C}$ with $l(c,{\mathcal C})\leq m-3$. Let ${\mathcal D}=(Q,\varphi|_Q)$ be the maximal sub-chord diagram of ${\mathcal C}$ such that $x,\varphi(x)\in Q$ and $l(c,{\mathcal D})=1$. Let $n$ be the number of chords of ${\mathcal D}$. 
Let $A$ (resp. $B$) be the point in $Q$ such that each components of ${\mathbb S}^1\setminus\{x,A\}$ (resp. ${\mathbb S}^1\setminus\{\varphi(x),B\}$) contains $n-1$ points of $Q$. Note that $n\geq m-(l(c,{\mathcal C})-1)\geq m-(m-3-1)=4$. Therefore we have that ${\mathcal D}$ has at least $4$ chords. Then there is a chord $d=(y,\varphi(y))$ of ${\mathcal D}$ such that $\{y,\varphi(y)\}$ and $\{x,\varphi(x),A,B\}$ are mutually disjoint. 
Suppose that $y$ and $\varphi(y)$ do not separate $x$ and $A$. In this case there must be a chord $e=(z,\varphi(z))$ of ${\mathcal D}$ such that $\{z,\varphi(z)\}$ and $\{x,\varphi(x),y,\varphi(y)\}$ are mutually disjoint and $z$ and $\varphi(z)$ do not separate $y$ and $\varphi(y)$. Then we have that the chords $c$, $d$ and $e$ form a sub-chord diagram of ${\mathcal C}$ that is equivalent to ${\mathcal C}_3$. This is a contradiction.
Suppose that $y$ and $\varphi(y)$ separate $x$ and $A$. Let $p$ and $q$ be a cutting pair of ${\mathcal C}\setminus d$. Then we have by Lemma \ref{lemma1} that $p$ and $q$ do not separate $y$ and $\varphi(y)$. Note that $p$ and $q$ is also a cutting pair of ${\mathcal D}\setminus d$ and they separate $x$ and $\varphi(x)$. Then we have that the component of ${\mathbb S}^1\setminus\{p,q\}$ that contains both $A$ and $B$ has more points of $Q\setminus\{y,\varphi(y)\}$ than the other. This is a contradiction. $\Box$

\vskip 3mm

Thus we have shown the following lemma.

\vskip 3mm

\begin{Lemma}\label{lemma4}
Let ${\mathcal C}=(P,\varphi)$ be an $m$-chord diagram on ${\mathbb S}^1$ that satisfies the condition $(\ast)$. 
Then ${\mathcal C}$ satisfies the following condition $(\star)$.

$(\star)$ $l(c,{\mathcal C})=m-2$ for every chord $c$ of ${\mathcal C}$.
\end{Lemma}

\vskip 3mm

\begin{Proposition}\label{proposition2}
Let ${\mathcal C}=(P,\varphi)$ be an $m$-chord diagram on ${\mathbb S}^1$ that satisfies the condition $(\star)$. 
If $m$ is even then $m$ is divisible by $4$ and ${\rm arc}({\mathcal C})=2$. 
If $m$ is odd then ${\mathcal C}$ is equivalent to ${\mathcal C}_m$.
\end{Proposition}

\vskip 3mm

\noindent{\bf Proof.}
Recall that $R_{2m}$ is a regular $(2m)$-gon inscribed in ${\mathbb S}^1$ and $v_{2m;1},\cdots,v_{2m;2m}$ are the vertices of $R_{2m}$ lying in this order. 
Let $G_{2m,m-2}$ be the graph whose vertices are $v_{2m;1},\cdots,v_{2m;2m}$ and whose edges are the chords $c(2m;i,m-2)$ joining the vertices $v_{2m;i}$ and $v_{2m;i+m-2}$ where $i\in\{1,\cdots,2m\}$. By calculating the greatest common divisor $(2m,m-2)=(2m-2(m-2),m-2)=(4,m-2)$ we have the isomorphism type of the graph $G_{2m,m-2}$ as follows.

(1) If $m$ is a multiple of $4$ then $(4,m-2)=2$ and therefore $G_{2m,m-2}$ is isomorphic to a disjoint union of two $m$-cycles. 

(2) If $m$ is congruent to $2$ modulo $4$ then $(4,m-2)=4$ and therefore $G_{2m,m-2}$ is isomorphic to a disjoint union of four $\frac{m}{2}$-cycles.

(3) If $m$ is odd then then $(4,m-2)=1$ and therefore $G_{2m,m-2}$ is isomorphic to a $2m$-cycle. 

\begin{figure}[htbp]
      \begin{center}
\scalebox{0.4}{\includegraphics*{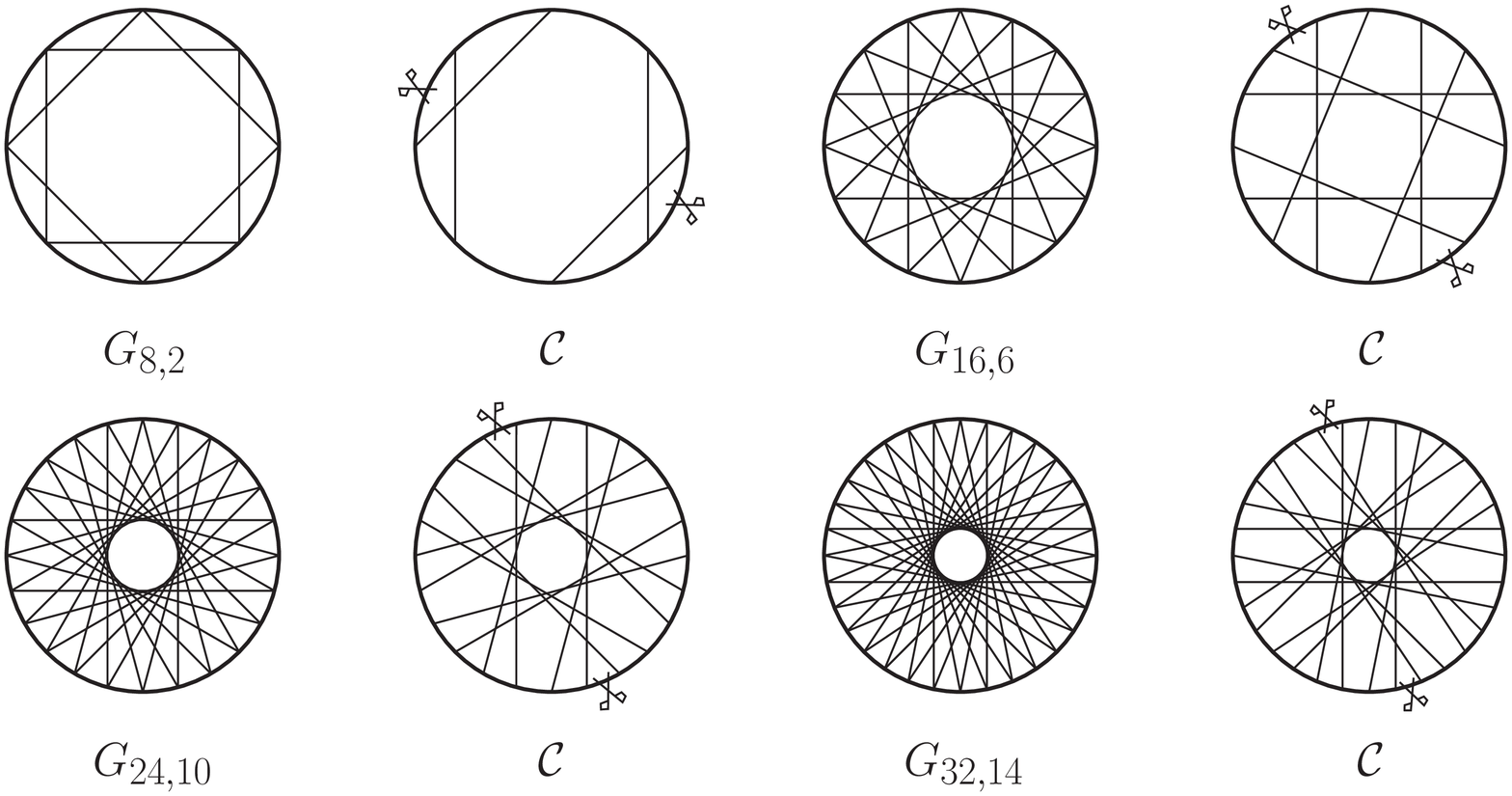}}
      \end{center}
   \caption{}
  \label{matching}
\end{figure}

\noindent
Note that in each case ${\mathcal C}$ must be a complete matching of the graph $G_{2m,m-2}$. In (1) we have up to symmetry that ${\mathcal C}$ is as illustrated in Figure \ref{matching}. Then we have that ${\rm arc}({\mathcal C})=2$. In (2) we have that $G_{2m,m-2}$ has no complete matchings because an $\frac{m}{2}$-cycle is an odd-cycle. In (3) we have that ${\mathcal C}$ is equivalent to ${\mathcal C}_m$. 
This completes the proof. $\Box$

%\vskip 3mm
% 

\vskip 3mm
\noindent{\bf Proof of Proposition \ref{proposition}.}
By Lemma \ref{lemma4} and Proposition \ref{proposition2} we have the result. $\Box$

\section{Examples of plane curves} 

Let ${\mathcal C}=(P,\varphi)$ be a chord diagram and $c=(x,\varphi(x))$ and $d=(y,\varphi(y))$ two chords of ${\mathcal C}$. We say that $c$ and $d$ are {\it parallel} if the pair of points $x$ and $\varphi(x)$ do not separate the pair of points $y$ and $\varphi(y)$. We say that two distinct points $x$ and $y$ in $P$ are {\it next to each other} if there is a component of ${\mathbb S}^1\setminus\{x,y\}$ that is disjoint from $P$. 
We say that $c$ and $d$ are {\it close to each other} if $x$ and $y$ are next to each other and $\varphi(x)$ and $\varphi(y)$ are next to each other, or $x$ and $\varphi(y)$ are next to each other and $\varphi(x)$ and $y$ are next to each other. 

\vskip 3mm

\begin{figure}[htbp]
      \begin{center}
\scalebox{0.47}{\includegraphics*{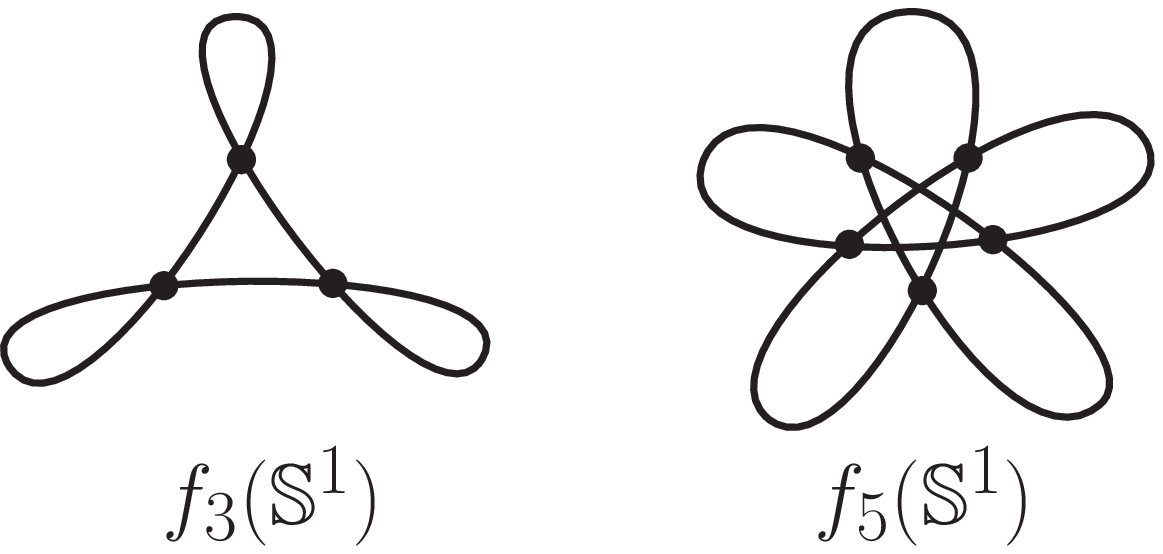}}
      \end{center}
   \caption{}
  \label{symmetric}
\end{figure}

\begin{figure}[htbp]
      \begin{center}
\scalebox{0.47}{\includegraphics*{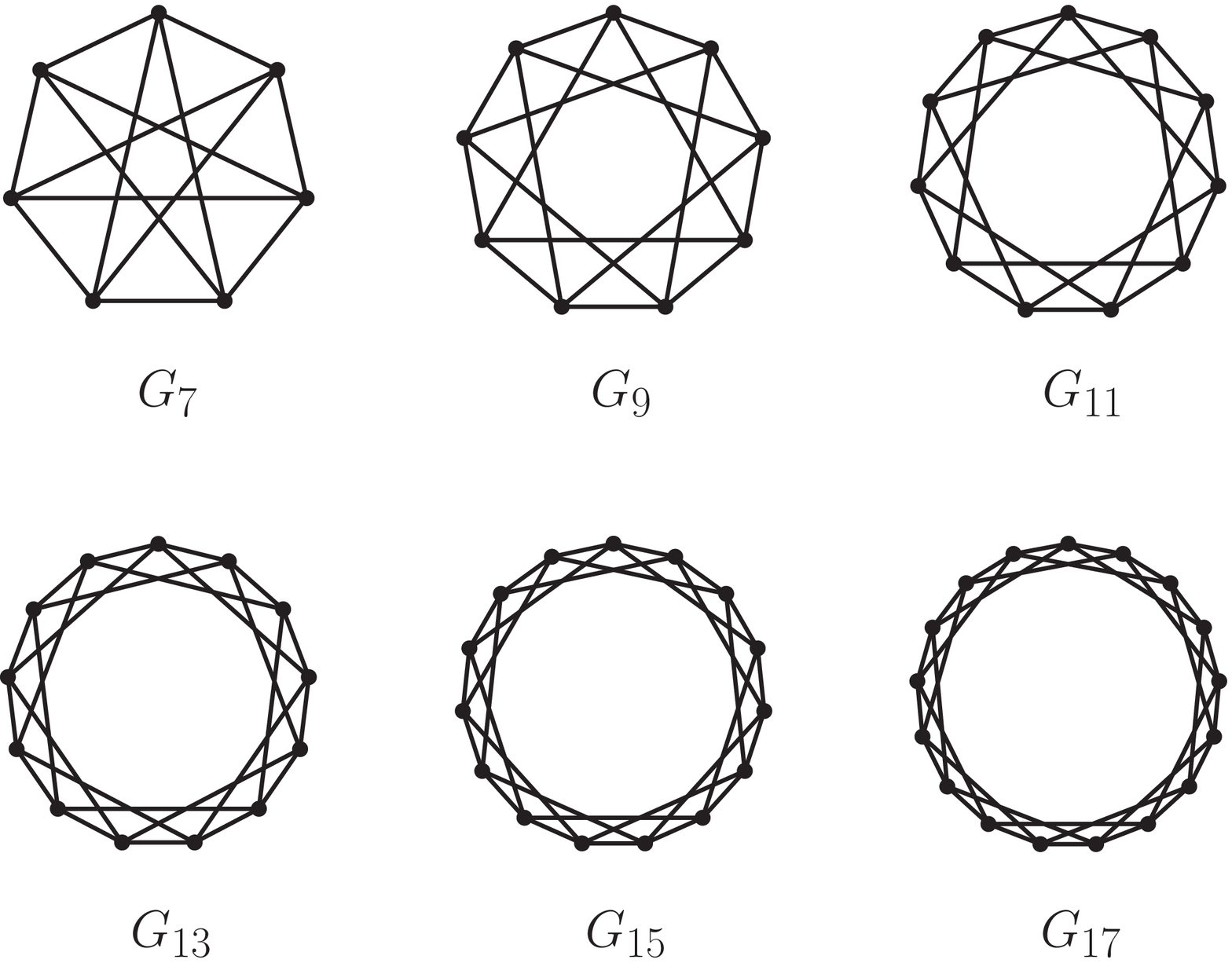}}
      \end{center}
   \caption{}
  \label{graph}
\end{figure}

\begin{figure}[htbp]
      \begin{center}
\scalebox{0.47}{\includegraphics*{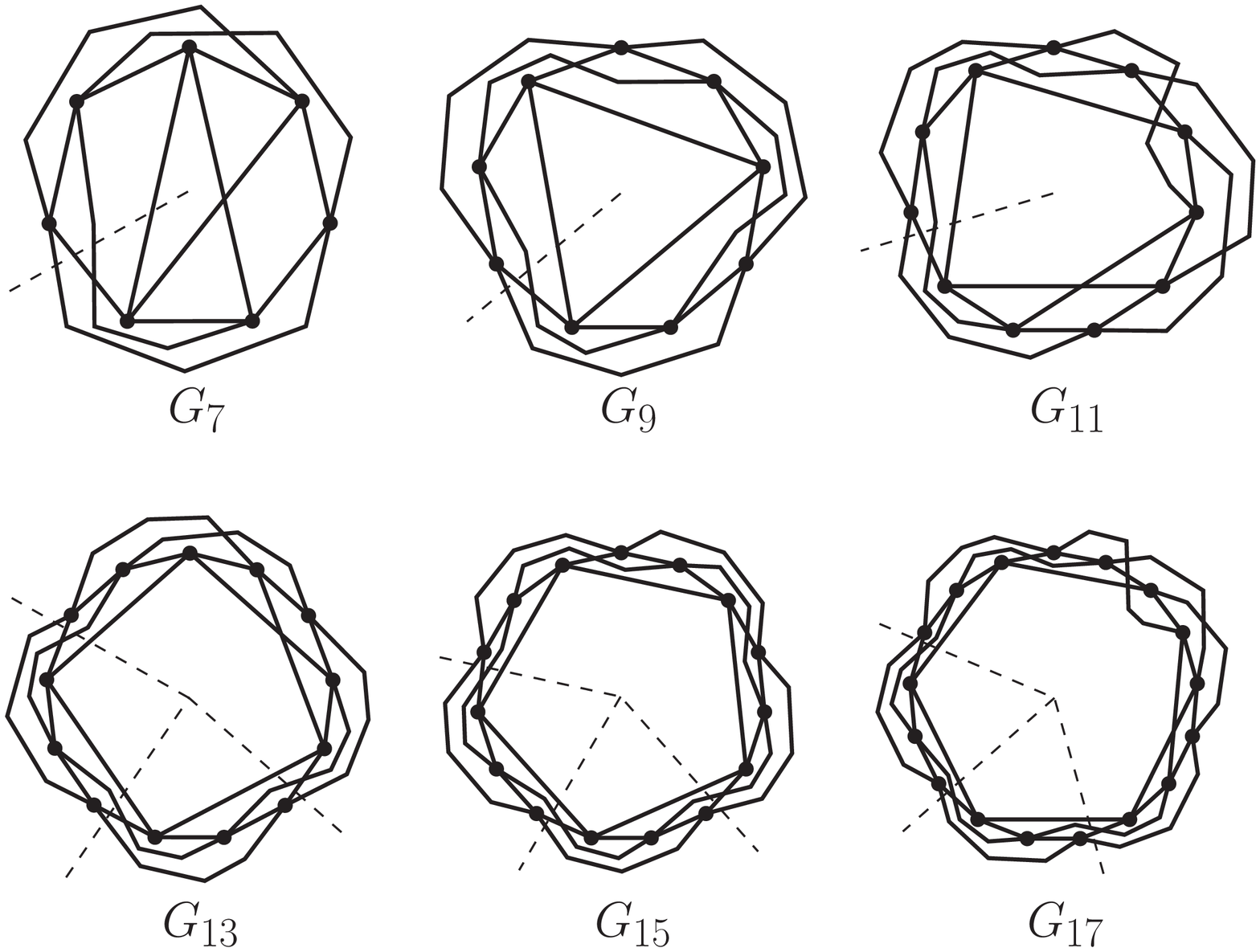}}
      \end{center}
   \caption{}
  \label{graph2}
\end{figure}

\begin{figure}[htbp]
      \begin{center}
\scalebox{0.47}{\includegraphics*{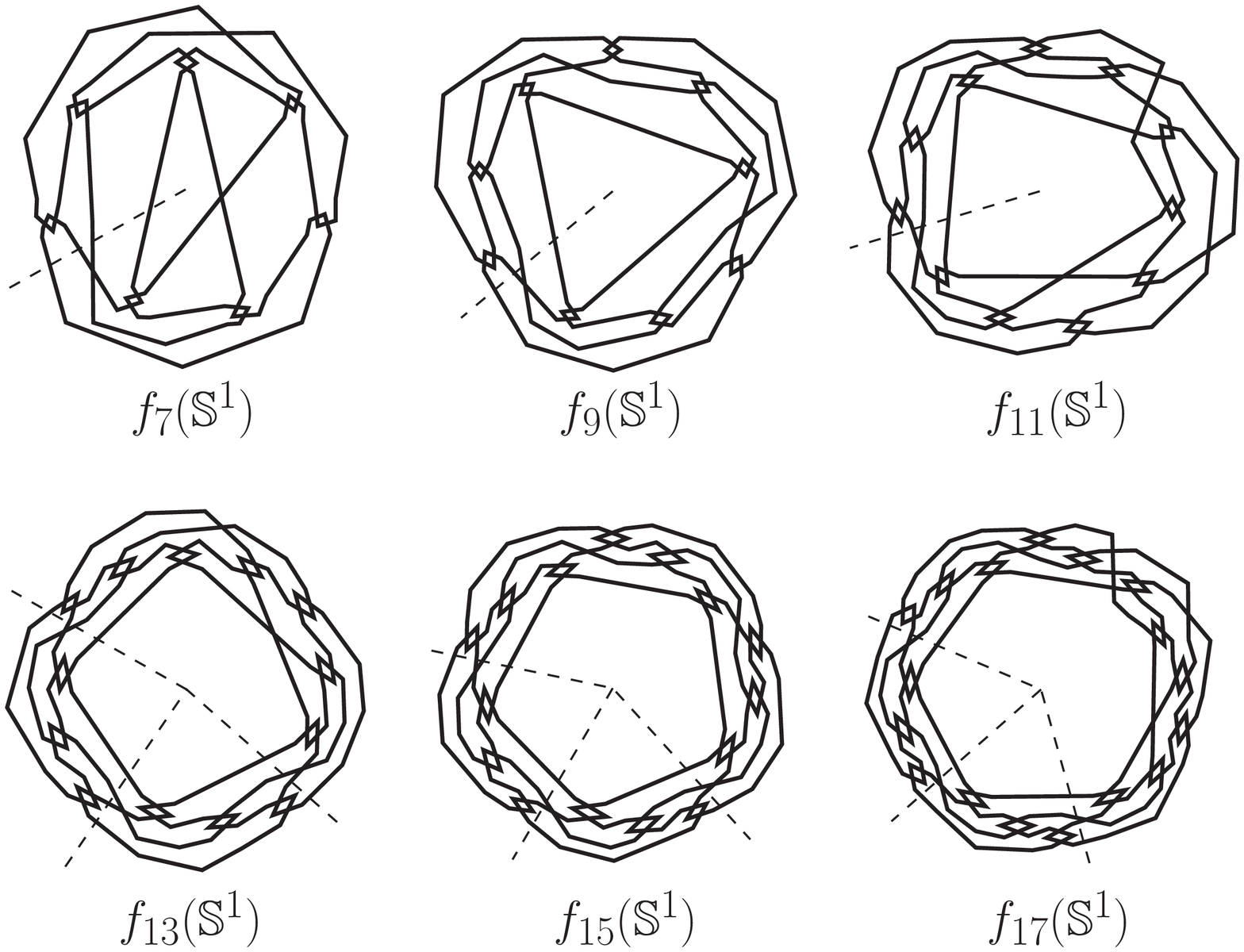}}
      \end{center}
   \caption{}
  \label{plane-curve}
\end{figure}

\noindent{\bf Proof of Proposition \ref{example}.}
The cases $n=1,2$ are shown in Figure \ref{symmetric}. We consider the case $n\geq3$. Let $G_{2n+1}={\mathbb S}^1/\sim_{{\mathcal C}_{2n+1}}$ be the $4$-regular graph obtained from ${\mathbb S}^1$ by identifying the end points of each chord of ${\mathcal C}_{2n+1}$. It is easy to observe that $G_{2n+1}$ is isomorphic to a graph obtained from a $(2n+1)$-cycle $\Gamma_{2n+1}$ on vertices $v_1,\cdots,v_{2n+1}$ lying in this order by adding edges joining $v_i$ and $v_{i+3}$ for each $i$ such that along the counterclockwise orientation of ${\mathbb S}^1$ the vertices of $G_{2n+1}$ appears $v_i,v_{i+1},v_{i+1-3},v_{i+1-3+1},v_{i+1-3+1-3},\cdots$. See Figure \ref{graph}. 
Then we deform them on ${\mathbb R}^2$ as illustrated in Figure \ref{graph2}. Note that they are classified into three types by $2n+1$ modulo $6$. 
Namely $G_{2n+1+6}$ is obtained from $G_{2n+1}$ by cutting open $G_{2n+1}$ along the dotted line and inserting two pieces of a pattern as illustrated in Figure \ref{graph2}.
We modify this $G_{2n+1}$ and have the image $f_n({\mathbb S}^1)$ as illustrated in Figure \ref{plane-curve}. Note that each vertex of $G_{2n+1}$ is replaced by two transversal double points. We call them a {\it twin pair}. The chords corresponding to them are also called a twin pair. 
By choosing any one of them for each twin pair we have a sub-chord diagram of ${\mathcal C}(f_n)$ that is equivalent to ${\mathcal C}_{2n+1}$ by the construction. Observe that each $f_n({\mathbb S}^1)$ is made of $(2n+1)$-times repetitions of \lq\lq one step forward and three steps back\rq\rq\  along the $(2n+1)$-cycle $\Gamma_{2n+1}$ and it totally goes around $\Gamma_{2n+1}$ twice. Here \lq\lq one step forward\rq\rq\ corresponds to an edge of $G_{2n+1}$ joining $v_i$ and $v_{i+1}$ and \lq\lq three steps back\rq\rq\  corresponds to an edge of $G_{2n+1}$ joining $v_{i+1}$ and $v_{i+1-3}$. 
It has no local double points and each double point comes from a part and another part that is one lap behind. 
Therefore we have that ${\rm arc}(f_n)=3$. 

Now we will check that no sub-chord diagram of ${\mathcal C}(f_n)$ is equivalent to ${\mathcal C}_{2m+1}$ for any $m<n$. 
Note that two chords in a twin pair are close to each other in ${\mathcal C}(f_n)$. Let ${\mathcal D}(f_n)$ be a sub-chord diagram of ${\mathcal C}(f_n)$ obtained from ${\mathcal C}(f_n)$ by deleting one of two chords for each twin pair in ${\mathcal C}(f_n)$. 
Since no two chords in ${\mathcal C}_{2m+1}$ are close to each other it is sufficient to check that no sub-chord diagram of ${\mathcal D}(f_n)$ is equivalent to ${\mathcal C}_{2m+1}$ for any $m<n$. 
Suppose that ${\mathcal E}$ is a sub-chord diagram of ${\mathcal D}(f_n)$ that is equivalent to ${\mathcal C}_{2m+1}$ for some $m<n$. Since no proper sub-chord diagram of ${\mathcal C}_{2n+1}$ is equivalent to ${\mathcal C}_{2m+1}$ we have that there is a chord $c$ of ${\mathcal E}$ that does not belong to any twin pair of ${\mathcal C}(f_n)$. Namely $c$ corresponds to a transversal double point of $f_{n}$ that comes from a double point of $G_{2n+1}\subset{\mathbb R}^2$ in Figure \ref{graph2}. 
Observe that for each chord $d=(x,\varphi(x))$ of ${\mathcal C}_{2m+1}$ there exist exactly two chords $g=(y,\varphi(y)$ and $h=(z,\varphi(z))$ of ${\mathcal C}_{2m+1}$ that are parallel to $d$ such that all of $y,\varphi(y),z$ and $\varphi(z)$ are contained in the same component of ${\mathbb S}^1\setminus\{x,\varphi(x)\}$. 
Therefore $c$ must have such two chords in ${\mathcal D}(f_n)$. By the \lq\lq one step forward and three steps back\rq\rq\  structure of $f_n({\mathbb S}^1)$ mentioned above the double points corresponding to such chords must lie in a small neighbourhood of the double point corresponding to $c$. Then we can check that there are no such two chords for $c$ in ${\mathcal D}(f_n)$ except the case $2n+1$ is congruent to $5$ modulo $6$ and $c$ is one of the three chords of ${\mathcal C}(f_n)$ that come from the three double points on the same edge of $G_{2n+1}\subset{\mathbb R}^2$ in Figure \ref{graph2}. 
For this exceptional case we further observe that the chords $g$ and $h$ above intersect unless $m=1$ and the end point of $g$ (resp. $h$) that is next to $x$ or $\varphi(x)$ in ${\mathcal C}_{2m+1}$ is not next to any end point of $h$ (resp. $g$) in ${\mathcal C}_{2m+1}$. 
However in this exceptional case we can check that there are no such two chords for $c$ in ${\mathcal D}(f_n)$. 
This is a contradiction. $\Box$

%\vskip 3mm
% 

%

%\vskip 3mm
% 

%

%\vskip 3mm
% 

%

%\vskip 3mm
% 

\vskip 3mm

\begin{Remark}\label{remark}
{\rm
It is easy to see that the graph $G_{2n+1}={\mathbb S}^1/\sim_{{\mathcal C}_{2n+1}}$ is a non-planar graph for $n\geq2$. Therefore we have that for $n\geq2$ there is no smooth immersion $f:{\mathbb S}^1\to {\mathbb R}^2$ that has only finitely many transversal double points whose associated chord diagram ${\mathcal C}(f)$ itself is equivalent to ${\mathcal C}_{2n+1}$. 
}
\end{Remark}

\vskip 3mm

\newpage
\ 

\section*{Acknowledgments} This work has been done while the author was visiting George Washington University. The author is grateful to the hospitality of the mathematical department of the George Washington University. The author is thankful to Professor J\'{o}zef Przytycki, Dr. Ken Shoda and Dr. Radmila Sazdanovi\'{c} for valuable discussions on mathematics.

{\normalsize
}

\end{document}